\documentclass{amsart}
\usepackage{amsfonts}
\usepackage{amsmath,amsthm}
\usepackage{amsfonts,amssymb}
\usepackage{mathrsfs}
\usepackage{enumerate}
\usepackage{verbatim}

\usepackage{multirow}
\usepackage{booktabs}

\usepackage[numbers,sort&compress]{natbib}

\newtheorem{thm}{Theorem}[section]
\newtheorem{cor}[thm]{Corollary}
\newtheorem{lem}[thm]{Lemma}

\newtheorem{exam}[thm]{Example}

\newtheorem{defn}[thm]{Definition}
\newtheorem{rem}{Remark}

\numberwithin{equation}{section}

\newcommand{\thmref}[1]{Theorem~\ref{#1}}

\newcommand{\lemref}[1]{Lemma~\ref{#1}}

\newcommand{\al}{\alpha}

\def\vz{\varepsilon}

\def\lz{\lambda}

\def\az{\alpha}

\def\({\Bigl(}
\def \){ \Bigr)}

\newcommand{\uu}{\mathfrak{u}}

\newcommand{\NN}{\mathbb{N}}

\newcommand{\bsgamma}{\boldsymbol{\gamma}}

\newcommand{\bsj}{\boldsymbol{j}}

 \def\NN{{\mathbb N}}

 \def\RR{{\mathbb R}}

\def\bsj{{\bf j}}

\def\wt{\widetilde}

\allowdisplaybreaks

\makeatletter
\@namedef{subjclassname@2020}{\textup{2020} Mathematics Subject Classification}
\makeatother

\begin{document}
\def\RR{\mathbb{R}}
\def\Exp{\text{Exp}}
\def\FF{\mathcal{F}_\al}

\title[] {Worst case tractability of weighted linear
tensor product problems under the absolute error criterion }

\author{Heping Wang} \address{ School of Mathematical Sciences, Capital Normal
University, Beijing 100048, China.}\email{wanghp@cnu.edu.cn}

\author{Yuheng Xi} \address{ School of of Mathematical Sciences, Capital Normal
University, Beijing 100048, China.} \email{2250502119@cnu.edu.cn.}


\begin{abstract}

This paper studies weighted linear tensor product problems in the
worst case setting. We consider algorithms that use finitely many
evaluations of arbitrary continuous linear functionals. We
investigate  algebraic and exponential tractability under the
absolute error criterion  in the case $\lambda_1>1$, where
$\lambda_1$ is the square of the univariate maximal singular
value. We obtain  necessary and sufficient conditions for various
notions of  algebraic and exponential tractability on weights
 and univariate singular values. Our approach is to convert  algebraic tractability of  the problems
for the absolute error criterion with $\lambda_1 > 1$ to
exponential tractability of the normalized problems. Specifically,
we solve Open Problem 27 as posed by Novak and Wo\'zniakowski
(2008).
\end{abstract}

\maketitle
\input amssym.def

\section{Introduction}\

In the field of information-based complexity, the computational
difficulty of a multivariate problem $S = \{S_d\}_{d\in \NN}$ is
measured by its information complexity $n(\vz, S_d)$, i.e., the
minimal number of information evaluations needed to approximate
$S_d$ within an error threshold $\vz > 0$, and tractability
studies how $n(\vz, S_d)$ depends on the dimension $d$ and on
$\vz^{-1}$; we refer to the monograph \cite{NW08, NW10, NW12} for
the general theory. Depending on how the error is scaled, one
distinguishes the absolute error criterion (ABS) and the
normalized error criterion (NOR), and depending on the way in
which $n(\vz, S_d)$ is allowed to grow, one distinguishes
algebraic (ALG) tractability, which regards $n(\vz, S_d)$ as a
function of $d$ and $\vz^{-1}$, and exponential (EXP)
tractability, in which $\vz^{-1}$ is replaced by $1 + \ln
\vz^{-1}$. Recently the study of ALG tractability and EXP
tractability has attracted much interest, and a great number of
interesting results
 have been obtained (see
\cite{CW1, DKPW14, DLPW11, GW11, IKPW16a, KW19, KPW, NW08, NW10,
NW12, PP14, PPW17, S13, SW15, W19, X15} and the references
therein).

This paper is devoted to linear weighted tensor product problems
with dimension-independent product weights $\bsgamma =
\{\gamma_j\}_{j\in \NN}$ in the worst case setting under the
absolute error criterion; the algorithms may use arbitrary
continuous linear functionals (the class $\Lambda^{\rm all}$).
Such a problem $S_{\bsgamma} = \{S_{d,\bsgamma}\}_{d\in \NN}$ is
determined by two non-increasing sequences: the sequence
$\{\lambda_j\}_{j\in \mathbb N}$, whose terms are the squares of
the ordered singular values of the univariate operator, and the
product weights $\{\gamma_j\}_{j\in \mathbb N}$, which moderate
the importance of the $j$-th univariate problem in the definition
of the tensor product for the multivariate case. Without loss of
generality we may assume that $\lz_2>0$ and that  all weights
$\gamma_j$ are positive. The unweighted case is obtained if we
take $\gamma_j = 1$ for all $j\in \mathbb N$.

Weighted linear tensor product problems of this kind arise
naturally in the study of multivariate approximation with tensor
product kernels and their tractability has been studied
extensively.

For unweighted linear tensor product problems and for algebraic
 tractability, algebraic polynomial tractability (ALG-PT) and algebraic strong
polynomial tractability (ALG-SPT) was characterized  by
Wo\'zniakowski \cite{W94};
 algebraic quasi-polynomial tractability (ALG-QPT) by Gnewuch and
Wo\'zniakowski \cite{GW11}; algebraic weak tractability (ALG-WT)
by Wo\'zniakowski \cite{W94} and by Papageorgiou and Petras
\cite{PP09};  the refined algebraic uniform weak tractability
(ALG-UWT) by Siedlecki \cite{S13};  the refined algebraic
$(s,t)$-weak tractability (ALG-$(s,t)$-WT) by Siedlecki and Weimar
\cite{SW15} and by Liu and Wang \cite{LW26}. Results on EXP
tractability were obtained in \cite{HKW19, PP14, PPW17}.

However, we have many negative results for algebraic polynomial
tractability and exponential tractability.  For example for
algebraic  tractability it is known from \cite[Theorem 5.5]{NW08}
that $\lambda_1 > 1$ leads to the curse of dimensionality, whereas
$\lambda_1 = 1$ leads to polynomial intractability. Under NOR,
which coincides with ABS when $\lambda_1 = 1$, we have polynomial
intractability. For exponential tractability it is known from
\cite{HKW19, PP14, PPW17} that exponential weak tractability does
not hold.

Now we consider ALG tractability and EXP tractability of  weighted
linear tensor product problems. For NOR algebraic tractability of
the problems was characterized in \cite[Section 5.3.4]{NW08}, and
exponential tractability of the problems with $\lz_1=1$ was
characterized by Kritzer, Pillichshammer and Wo\'zniakowski
\cite{KPW} except for EXP-$(s,t)$-WT with $\max(s,t) < 1$ and
EXP-UWT, which were settled by Liu, Wang and Wang \cite{LWW27}.

For ABS with $\lambda_1 \le 1$ algebraic tractability of the
problems was essentially characterized in \cite[Sections 5.2 and
5.3]{NW08}. For ABS with $\lambda_1 > 1$ the situation  is much
less understood. In this case the initial error $\lambda_1^{d/2}$
grows exponentially with the dimension $d$, ABS differs
essentially from NOR. Novak and Wo\'zniakowski posed Open Problem
27 in \cite{NW08}:  find necessary and sufficient conditions for
algebraic polynomial and weak tractability of $S_{\bsgamma}$ for
ABS in the worst case setting.

In this paper for  weighted linear tensor product problems we give
necessary and sufficient conditions for various notions of ALG
tractability and EXP tractability for ABS with $\lz_1>1$.
Specifically, we solve Open Problem 27 completely for
dimension-independent product weights under ABS with $\lz_1>1$. We
show that $S_{\bsgamma}$ is  {\it never} ALG-SPT for arbitrary
weights, and that
     $S_{\bsgamma}$ is ALG-PT if and only if
\begin{equation*}
\lim_{j \to \infty} \gamma_j = 0 \quad \mbox{and} \quad
\limsup_{j+k \to \infty} \frac{k \ln j}{\ln \ln (\gamma_k^{-1}+
\lambda_j^{-1})} <\infty.
\end{equation*} We also show that $S_{\bsgamma}$ is ALG-WT if and only
if
\begin{equation*}
\lim_{j \to \infty} \gamma_j = 0 \quad \mbox{and} \quad \lim_{j
\to \infty} \frac{\ln \lambda_j^{-1}}{\ln j} = \infty.
\end{equation*}
See Definition \ref{def:M},   Theorems \ref{thm:main} and
\ref{cor:M} for details.

Our approach is to convert  ALG tractability of  the problems for
ABS with $\lambda_1 > 1$ to  EXP tractability of the normalized
problems $\wt{S}_{\bsgamma}$. The critical thought comes from
Remark 1 of \cite{LW26}, where it is observed that, for unweighted
linear tensor product problems and for ABS with $\lambda_1
> 1$, ALG-$(s,t)$-WT is equivalent to EXP-$(t,t)$-WT with $t>1$ and to
EXP-$(s,t)$-WT with $s > t$. The key is to find  the inequality
about the information complexity with different dimension, see
\lemref{lem:padding}.

The paper is organized as follows. In Sections 2.1 and 2.2 we
introduce  weighted tensor product problems and various notions of
algebraic and exponential tractability. Section 2.3 collects known
results about ALG tractability in the case $\lambda_1 \le 1$ and
about EXP tractability for the normalized problem. Section 2.4
gives  main results.  Section 3 contains  proofs of the main
results.

\section{Preliminary}

\subsection{Weighted linear tensor product problems}\

Let $S_1 : H_1 \to G_1$ be a compact linear operator between
separable Hilbert spaces.  Then $W_1 := S_1^* S_1: H_1 \to H_1$ is
a compact self-adjoint non-negative operator, where $T^*$ denotes
the adjoint
 of a linear operator $T$. Let $\{(\lz_j, e_j)\}_{j\in
\mathbb N}$ denote the eigenpairs of $W_1$, that is,
$$W_1e_j=\lz_je_j,\ e_j\in H_1,\ \ \lz_1\ge \lz_2\ge\dots\ge 0,$$
and $\{e_j\}_{j\in\mathbb N}$ is an orthonormal basis in $H_1$.
Without loss of generality, we assume that $\lambda_2>0$ in the
sequel to exclude the trivial rank-one case.  By compactness of
$S_1$, we have $\lim\limits_{j\to\infty}\lambda_j=0$.

For a real number $\gamma > 0$, let $H_{1,\gamma}$ be the Hilbert
space $H_1$ endowed with the weighted inner product
\begin{equation*}
\langle f, g \rangle_{H_{1,\gamma}} := \langle f_0, g_0
\rangle_{H_1} + \gamma^{-1} \langle f_1, g_1 \rangle_{H_1},
\end{equation*}
where $f = f_0 + f_1$ is the decomposition along $H_{1,0} := {\rm
span}\, \{e_1\}$ and its orthogonal complement. The weighted
operator $S_{1,\gamma}$ is defined as the operator $S_1$ from
$H_{1,\gamma}$ to $G_1$.

Clearly, $S_{1,\gamma}:H_{1,\gamma} \to G_1$ and $W_{1,\gamma} =
S_{1,\gamma}^* S_{1,\gamma}$ are compact linear operators, and the
eigenpairs of $W_{1,\gamma}$ are $\{(\lz_{\gamma,j}, e_j)\}_{j\in
\mathbb N}$, where
\begin{equation*}
\lambda_{\gamma,j}:= \left\{
\begin{array}{ll}
\lambda_1\,, \quad & j = 1\; ,\\
\gamma \lambda_j\,, \quad & j \ge 2\;.
\end{array}\right.
\end{equation*}

Now we introduce the {\it dimension-independent product weights}
$\bsgamma = \{\gamma_j\}_{j\in \NN}$. Without loss of generality,
we assume
\begin{equation*}
1 \ge \gamma_1 \ge \gamma_2 \ge \dots > 0.
\end{equation*}
  For $d\in \NN$, the $d$ times
tensor product spaces $H_{d,\bsgamma}$ and $G_d$ are defined by
$$H_{d,\bsgamma}=H_{1,\gamma_1}\otimes
H_{1,\gamma_2}\otimes\dots\otimes H_{1,\gamma_d},\ \
G_d=G_1\otimes \dots\otimes G_1,$$ respectively, the tensor
product operator $S_{d,\bsgamma}: H_{d,\bsgamma}\to G_d$ is
represented by
$$S_{d,\bsgamma}=S_{1,\gamma_1}\otimes\dots\otimes
S_{1,\gamma_d}.$$ The weighted linear tensor product problem is
referred to $S_{\bsgamma}=\{S_{d,\bsgamma}\}_{d\in\mathbb N}$.

If $\gamma_j\equiv1$ for all $j\in\NN$, then the weighted tensor
product space $H_{d,\bsgamma}$  recedes to the tensor product
space $H_d=H_1\otimes \cdots \otimes H_1$;  the weighted linear
tensor product operator $S_{d,\bsgamma}$ reduces to the unweighted
operator $S_d=S_1\otimes \cdots \otimes S_1$;  the weighted linear
tensor product problem recedes to the (unweighted) linear tensor
product problem.

 Clearly,
$W_{d,\bsgamma}=S_{d,\bsgamma}^*S_{d,\bsgamma}: H_{d,\bsgamma}\to
H_{d,\bsgamma} $ is a compact linear operator and the eigenpairs
of $W_{d,\bsgamma}$ are
$\{(\lz_{d,\bsgamma,\bsj},e_{\bsj})\}_{\bsj\in\mathbb N^d}$, where
$$\lz_{d,\bsgamma,\bsj}=\prod_{k=1}^d \lz_{\gamma_k,j_k},  \ \ \  \ e_{\bsj}=e_{j_1}\otimes\cdots \otimes e_{j_d}.$$

We write $[d]:=\{1,2,\dots,d\},$ $\uu(\bsj):=\{k \in [d]: j_k>1\}$
for $ \bsj= (j_1, \dots, j_d) \in \NN^d$, and
\begin{equation*}
\lambda_{d,\bsj}:=\prod_{k=1}^{d} \lambda_{j_k},\ \
\gamma_{\uu(\bsj)} :=\prod_{k \in \uu(\bsj)} \gamma_{k}.
\end{equation*}
 For  $ \bsj \in \NN^d$, the eigenvalues of
$W_{d,\bsgamma}$ can also be written as
\begin{equation}\label{eq:eig2}
\lambda_{d,\bsgamma,\bsj} =\prod_{k=1}^{d} \lambda_{\gamma_k, j_k}
=\gamma_{\uu(\bsj)} \prod_{k=1}^{d}
\lambda_{j_k}=\gamma_{\uu(\bsj)} \lambda_{d,\bsj}.
\end{equation}

We now consider  the linear weighted tensor product problem
$S_{\bsgamma}=\{S_{d,\bsgamma}\}_{d\in\NN}$. We want to
approximate $S_{d,\bsgamma}$ by algorithms $A_{d,n}(f)$ of the
form
\begin{equation}\label{1.2}
A_{d,n}(f)=\phi_{d,n}(L_1(f),L_2(f),\dots,L_n(f)),
\end{equation}
where  $L_1,\dots ,L_n$ belong to $\Lambda_d^{\rm
all}=H_{d,\bsgamma}^*$, the space of all continuous linear
functionals on $H_{d,\bsgamma}$, and $\phi_{d,n}:\RR^n \to G_d$ is
an arbitrary  mapping.  The worst case error of an algorithm
$A_{d,n}$ is defined as
$$ e^{\rm wor}(A_{d,n})=\sup_{\substack{f \in H_{d,\bsgamma} \\ \|f\|_{H_{d,\bsgamma}} \le 1}}\|S_{d,\bsgamma}(f)-A_{d,n}(f)\|_{G_d}.$$
The $n$-th minimal worst case error is defined by
$$e^{\rm wor}(n,S_{d,\bsgamma})= \inf_{A_{d,n}}e^{\rm wor}(A_{d,n}),$$ where the infimum is taken over
all algorithms of the form \eqref{1.2}. Let the sequence
$\{\lz_{d,j}\}_{j\in\mathbb N}$ be the non-increasing
rearrangement of $\{\lambda_{d,\bsgamma,\bsj}\}_{j\in{\mathbb
N}^d}$. It is well known that $$e^{\rm
wor}(n,S_{d,\bsgamma})=\sqrt{\lz_{d,n+1}}.$$
 When $n=0$, the initial
error $e_0({S}_{d,\bsgamma})$ is denoted by
$$ e_0({S}_{d,\bsgamma})=\sup_{\substack{f \in H_{d,\bsgamma} \\ \|f\|_{H_{d,\bsgamma}} \le 1}}\|S_{d,\bsgamma}(f)\|_{G_d}. $$
For $\vz\in(0,1)$, the information complexity for ABS or NOR is
defined by
\begin{equation*}
n^X (\vz,S_{d,\bsgamma}) :=\min\{n \in \NN: e^{\rm wor}
(n,S_{d,\bsgamma})\le \vz\,{\rm CRI}_d\},
\end{equation*}
where
\begin{equation*}
{\rm CRI}_d:=\left\{
\begin{array}{ll}
1\,, & X={\rm ABS}\; ,\\
e_0({S}_{d,\bsgamma})\,, & X={\rm NOR}\;.
\end{array}\right.
\end{equation*}
In our setting we have
\begin{equation*}
e_0({S}_{d,\bsgamma})=\|S_{d,\bsgamma}\| =\|W_{d,\bsgamma}\|^{1/2}
=\max_{\bsj \in \NN^d}
\lambda_{d,\bsgamma,\bsj}^{1/2}=\lambda_1^{d/2},
\end{equation*}
and
\begin{equation*}
n^{\rm ABS}(\vz, S_{d,\bsgamma}) = \#\{ \bsj:=(j_1,\dots,j_d) \in
\NN^d: \lambda_{d,\bsgamma,\bsj}=\prod_{k=1}^{d}
\lambda_{\gamma_k, j_k}
> \vz^2 \}.
\end{equation*}
The normalized case $\lambda_1 = 1$, in which the two criteria
coincide, has been completely solved (see \cite{NW08}, \S 5.3.4),
hence throughout the paper we assume
\begin{equation*}
\lambda_1 > 1.
\end{equation*}
In this case the initial error grows exponentially with $d$,
and the absolute error criterion (ABS) differs essentially from
the normalized error criterion (NOR).

For brevity, in this paper we write
\begin{equation*}
n(\vz, d) := n^{\rm ABS}(\vz, S_{d,\bsgamma}).
\end{equation*}
Evidently, $n(\vz, d)$ is non-increasing in $\vz$. We are
interested in the behavior of $n(\vz,d)$ as $d \to \infty$ and
$\vz \to 0^+$.

We next introduce the normalized problem of a weighted linear
tensor problem. This notation will be used throughout the paper.

\begin{defn}\label{def:normalized}
Given a weighted linear tensor problem $S_{\bsgamma} =
\{S_{d,\bsgamma}\}$ with $\lambda_1 > 1$, its normalized problem
$\wt{S}_{\bsgamma} = \{\wt{S}_{d,\bsgamma}\}$ is given by
\begin{equation*} \wt{S}_{d,\bsgamma}:=
\lambda_1^{-d/2}{S}_{d,\bsgamma}.
\end{equation*}\end{defn}

It follows that  $\wt{W}_{d,\bsgamma} = \wt{S}_{d,\bsgamma}^*
\wt{S}_{d,\bsgamma}= \lz_1^{-d}{S}_{d,\bsgamma}^*
{S}_{d,\bsgamma}=\lz_1^{-d}{W}_{d,\bsgamma}$ has eigenvalues
$\wt{\lambda}_{d,\bsgamma,\bsj} =\gamma_{\uu(\bsj)}
\prod_{k=1}^{d}
\wt{\lambda}_{j_k}=\lambda_1^{-d}\lambda_{d,\bsgamma,\bsj},$ where
\begin{equation*}
\wt{\lambda}_j:= \lambda_j/\lambda_1, \quad 1 = \wt{\lambda}_1
\ge \wt{\lambda}_2 \ge \dots \ge 0, \quad \wt{\lambda}_2 > 0,
\end{equation*}
Hence, the initial error of the normalized problem is
$e_0(\wt{S}_{d,\bsgamma})=\|\wt{S}_{d,\bsgamma}\| =
\wt{\lambda}_1^{d/2} = 1$, so ABS and NOR coincide. We note that
\begin{align}
&\quad \ \wt{n}(\vz, d):=n^{\rm NOR} (\vz, {S}_{d,\bsgamma})
=n^{\rm ABS} (\vz, \wt{S}_{d,\bsgamma})=n^{\rm NOR} (\vz,
\wt{S}_{d,\bsgamma})\notag\\
&= \#\{j \in
\NN^d:\wt{\lambda}_{d,\bsgamma,\bsj}>\vz^2\}=\#\{j \in
\NN^d:{\lambda}_{d,\bsgamma,\bsj}>\lambda_1^d\vz^2\}
=n(\lambda_1^{d/2}\vz,d).
\label{eq:conv}
\end{align}

\subsection{Notions of tractability}\

 In this subsection we briefly recall the various
tractability notions for weighted linear tensor product problems
under ABS in the worst case setting. First we introduce all
notions of algebraic tractability.

\begin{defn}\label{def:M}
The problem $S_{\bsgamma} = \{S_{d,\gamma}\}$ is said to be
\begin{itemize}
\item {\it algebraically strongly polynomially tractable}
(ALG-SPT) if there exist $C, p \ge 0$ such that
\begin{equation*}
n(\vz, d)\le C \vz^{-p}, \quad \forall d \in \NN, \forall \vz \in
(0,1);
\end{equation*}
\item {\it algebraically polynomially tractable} (ALG-PT) if
there exist $C, q, p \ge 0$ such that
\begin{equation*}
n(\vz, d) \le C d^q \vz^{-p}, \quad \forall d \in \NN, \forall \vz
\in (0,1);
\end{equation*}
\item {\it algebraically quasi-polynomially tractable} (ALG-QPT)
if there exist $C, t \ge 0$ such that
\begin{equation*}
n(\vz, d) \le C \exp(t (1 + \ln \vz^{-1}) (1 + \ln d)) , \quad
\forall d \in \NN, \forall \vz \in (0,1);
\end{equation*}
\item {\it algebraically $(s,t)$-weakly tractable}
(ALG-$(s,t)$-WT) for fixed $s, t > 0$ if
\begin{equation*}
\lim_{\vz^{-1} + d \to \infty} \frac{\ln n(\vz, d)}{d^t +
\vz^{-s}} = 0;
\end{equation*}
If $s = t = 1$, then we simply say {\it algebraically weakly
tractable} (ALG-WT); \item {\it algebraically uniformly weakly
tractable} (ALG-UWT) if it is ALG-$(s,t)$-WT for all $s, t > 0$.
\end{itemize}
\end{defn}

If the problem $S_{\bsgamma}$ is not ALG-WT, then $S_{\bsgamma}$
is called ALG intractable. We say that the problem  $S_{\bsgamma}$
suffers from {the curse of dimensionality} if there exist positive
numbers $C, \,\vz _0,\, a $ such that for all $0<\vz\leq \vz _{0}$
and infinitely many $d\in \mathbb N$,
\begin{equation*}
n(\vz ,d)\geq C(1+a )^{d}.
\end{equation*}
Clearly, if $S_{\bsgamma}$ suffers from {the curse of
dimensionality}, then $S_{\bsgamma}$ is  ALG intractable.

If the $n$-th minimal error is exponentially convergent, then we
should study tractability with $\vz^{-1}$ being replaced by $(1
+\ln\vz^{-1})$), which is called exponential tractability. We now
give the above notions of exponential tractability in detail.

\begin{defn}
The problem $S_{\bsgamma}$ is said to be
\begin{itemize}
\item {\it exponentially strongly polynomially tractable}
(EXP-SPT) if there exist $C, p \ge 0$ such that
\begin{equation*}
n(\vz, d) \le C (1 + \ln \vz^{-1})^p, \quad \forall d \in \NN,
\forall \vz \in (0,1);
\end{equation*}
\item {\it exponentially polynomially tractable} (EXP-PT) if
there exist $C, q, p \ge 0$ such that
\begin{equation*}
n(\vz, d) \le C d^q (1 + \ln \vz^{-1})^p, \quad \forall d \in \NN,
\forall \vz \in (0,1);
\end{equation*}
\item {\it exponentially quasi-polynomially tractable} (EXP-QPT)
if there exist $C, t \ge 0$ such that
\begin{equation*}
n(\vz, d) \le C \exp(t (1 + \ln d) (1 + \ln(1 + \ln \vz^{-1}))),
\quad \forall d \in \NN, \forall \vz \in (0,1);
\end{equation*}
\item {\it exponentially $(s,t)$-weakly tractable}
(EXP-$(s,t)$-WT) for fixed $s, t > 0$ if
\begin{equation*}
\lim_{d+\vz^{-1} \to \infty} \frac{\ln n(\vz, d)}{d^t+(1+\ln
\vz^{-1})^s}=0;
\end{equation*}
If $s = t = 1$, then we say {\it exponentially weakly tractable}
(EXP-WT); \item {\it exponentially uniformly weakly tractable}
(EXP-UWT) if it is EXP-$(s,t)$-WT for all $s, t > 0$.
\end{itemize}
\end{defn}

Clearly,  for $0<s_1<s, \ 0<t_1<t$,   $(s_1,t_1)$-WT
$\Longrightarrow$ $(s,t)$-WT.  UWT holds if and only if $(t,t)$-WT
holds for all $t\in(0,1)$.  For ALG and EXP tractability we also
have

\begin{center}SPT $\Longrightarrow$ PT  $\Longrightarrow$  QPT
$\Longrightarrow$ UWT  $\Longrightarrow$ $(s,t)$-WT with
$s,t>0$.\end{center}

\subsection{Known Results on Tractability}\

In this subsection we introduce previous results on ALG and EXP
tractabilities and present main results. To shorten the notation
we write ``iff'' instead of ``if and only if''.

For a positive sequence ${\bf x} = \{x_j\}$ we denote its
sum-exponent by
\begin{equation*}
p_{\bf x} := \inf\{\tau > 0: \sum_{j=1}^{\infty} x_j^{\tau} <
\infty\},
\end{equation*}with the convention that $ \inf \emptyset :=
\infty$. Then $p_{\bf x }<\infty$ iff there exists two positive
constants $C$ and $r$ such that $x_j\le Cj^{-r}$ for all $j\in
\mathbb N$ iff $${\bf x}_*=\liminf\limits_ {j\to\infty}\frac{\ln
x_j^{-1}}{\ln j}>0.$$ In this case, $p_{\bf x}=({\bf x}_*)^{-1}$.

Under ABS  for the ALG tractability of the weighted  problem
$S_{\bsgamma}$  in the cases $\lambda_1 < 1$ and $\lambda_1 = 1$,
necessary and sufficient conditions for ALG-SPT and ALG-PT were
given in \cite{NW08}, for ALG-QPT in \cite{NW12} and \cite{GW11},
for ALG-$(s,t)$-WT in \cite{NW08}, \cite{PP09}, and \cite{SW15},
and for ALG-UWT in \cite{S13}. We remark that if $\lz_1\le 1$,
then $S_{\bsgamma}$
      does not suffers from the curse of dimensionality iff there exists
      a $j_0\in \mathbb N$ such that $\gamma_{j_0}\lz_2<1$.
  See the following ALG
tractability results of $S_{\bsgamma}$ for ABS.\vspace{2mm}

     Let $\lambda_1 < 1$.

   $\bullet$ For arbitrary weights $\{\gamma_j\}$,
 ALG-SPT holds iff ALG-PT holds iff ALG-QPT holds iff
$
p_{\lambda} < \infty ;
$

$\bullet$ ALG-UWT holds iff
\begin{equation*}
\lim_{j \to \infty} \frac{\ln \lambda_j^{-1}}{\ln \ln j} = \infty;
\end{equation*}

$\bullet$ For fixed $s, t > 0$,  ALG-$(s,t)$-WT holds iff
\begin{equation*}
\lim_{j \to \infty} \lambda_j (\ln j)^{2/s} =0.
\end{equation*}

\vspace{2mm}

     Let $\lambda_1 = 1$.

     $\bullet$ ALG-SPT holds iff ALG-PT holds iff \begin{equation*}
p_{\lambda} < \infty \quad \mbox{and} \quad p_{\gamma} < \infty;
\end{equation*}

      $\bullet$ ALG-QPT holds iff $p_\lz<\infty$ and there exists
      a $j_0\in \mathbb N$ such that $\gamma_{j_0}\lz_2<1$;

$\bullet$ ALG-UWT holds iff \begin{equation*} \lim_{j \to \infty}
\frac{\ln \lambda_j^{-1}}{\ln \ln j} = \infty,
\end{equation*} and there exists
      a $j_0\in \mathbb N$ such that $\gamma_{j_0}\lz_2<1$;

$\bullet$ ALG-$(s,t)$-WT holds with $t\in (0,1]$ and $s>0$ iff
\begin{equation*} \lim_{j \to \infty} \lambda_j (\ln j)^{2/s} =0,
\end{equation*}and there exists
      a $j_0\in \mathbb N$ such that $\gamma_{j_0}\lz_2<1$;

$\bullet$ ALG-$(s,t)$-WT holds with $t>1$ and $s>0$  iff
\begin{equation*} \lim_{j \to \infty} \lambda_j (\ln j)^{2/s} =0.
\end{equation*}

\vspace{1mm}

With $\vz \in (0,1)$  define  the functions
\begin{equation*}
j(\vz):=\max\{j \in \NN:\lambda_j>\lambda_1 \vz^2\}, \quad
d(\vz):=\sup\{d \in \NN:\gamma_d>\vz^2\},
\end{equation*}
with the convention $\sup \emptyset = 0$. Then $d(\vz)$ is
well-defined, and if $\gamma_j \to 0$ and all weights are
positive, then $d(\vz)$ is finite and tends to $\infty$ as $\vz
\to 0$; note that $j(\vz) \ge 1$ always holds since $\lambda_1 >
\lambda_1 \vz^2$.

 Under ABS  for the EXP tractability of the normalized weighted  problem
$\wt{S}_{\bsgamma}$, necessary and sufficient conditions for
EXP-SPT,  EXP-PT, EXP-QPT, and EXP-$(s,t)$-WT with $\max(s,t)\ge1$
were given in \cite{KPW}, and for EXP-$(s,t)$-WT with
$\max(s,t)<1$  and EXP-UWT in \cite{LWW27}.
  See the following EXP
tractability results of $\wt{S}_{\bsgamma}$ under ABS.\vspace{2mm}

$\bullet$  EXP-SPT holds iff EXP-PT holds iff
\begin{equation}\label{eq:exppt} \lim\limits_{j \to
\infty} \gamma_j = 0\ \ {\rm and} \quad \limsup_{\vz \to 0}
\frac{d(\vz) \ln j(\vz)}{\ln \ln \vz^{-1}} < \infty.
\end{equation}

 \vskip 1mm

 $\bullet$  EXP-QPT holds
 iff  $\lim\limits_{j \to \infty}
\gamma_j = 0$ and
$$
\limsup_{\varepsilon \rightarrow 0} \frac{d(\varepsilon) \ln
j(\varepsilon)}{[\ln d(\varepsilon)] \ln\ln
\frac{1}{\varepsilon}}< \infty.
$$

 \vskip 1mm

  $\bullet$   EXP-UWT holds iff
\begin{equation}\label{eq:expuwt}
\lim_{j \to \infty} \frac{\ln \ln \wt{\lambda}_j^{-1}}{\ln \ln j}
= \infty \quad \mbox{and} \quad \lim_{j \to \infty} \frac{\ln \ln
\gamma_j^{-1}}{\ln j}=\infty.
\end{equation}\vskip 1mm

 $\bullet$  Let $s=1$ and $t=1$. EXP-WT holds  iff
 \begin{equation}\label{1.8}
\lim_{j \rightarrow \infty} \gamma_j=0\ \ \mbox{ and }\ \ \lim_{j
\rightarrow \infty} \frac{\ln \wt{\lambda}_j^{-1}}{\ln j}=\infty.
\end{equation}
 \vskip 2mm

 $\bullet$  Let $s=1$  and $t<1$.
EXP-$(1,t)$-WT holds  iff
\begin{equation}\label{1.9}
\lim_{j \rightarrow \infty} \frac{\ln {\gamma_j}^{-1}}{\ln
j}=\infty\ \ {\rm and }\ \ \lim_{j \rightarrow \infty} \frac{\ln
\wt{\lambda}_j^{-1}}{\ln j}=\infty.
\end{equation}
 \vskip 2mm

 $\bullet$  Let $s=1$ and $t>1$.
EXP-$(s,t)$-WT holds  iff
\begin{equation}\label{eq:exp1gt}
\mbox{$\gamma_j$'s are arbitrary}\ \ \ \mbox{and}\ \ \
 \lim_{j \rightarrow \infty} \frac{\ln \wt{\lambda}_j^{-1}}{\ln
   j}=\infty.
\end{equation}
 \vskip 2mm

 $\bullet$  Let $s>1$, $t\le1$ and $\lambda_2<1$.
EXP-$(s,t)$-WT holds  iff
\begin{equation*}\label{eq:expsgt1a}
\mbox{$\gamma_j$'s are arbitrary}\ \ \ \mbox{and}\ \ \ \lim_{j
\rightarrow \infty}\frac{(\ln\wt{\lambda}_j^{-1})^s}{\ln
  j}=\infty.
\end{equation*}
 \vskip 2mm

 $\bullet$  Let $s>1$, $t\le1$ and $\lambda_2=1$.
EXP-$(s,t)$-WT holds  iff
\begin{equation*}\label{eq:expsgt1b}
\exists\,p\in\NN \ \ \mbox{with}\ \ \gamma_p<1\ \ \ \mbox{and}\ \
\ \lim_{j \rightarrow \infty}\frac{(\ln\wt{\lambda}_j^{-1})^s}{\ln
  j}=\infty.
\end{equation*}
 \vskip 2mm

 $\bullet$ Let $s>1$ and $t>1$.
EXP-$(s,t)$-WT holds  iff
\begin{equation}\label{eq:expsgt}
\mbox{$\gamma_j$'s are arbitrary}\ \ \ \mbox{and}\ \ \ \lim_{j
\rightarrow \infty}\frac{(\ln\wt{\lambda}_j^{-1})^s}{\ln
  j}=\infty.
\end{equation}
 \vskip 2mm

 $\bullet$  Let $s<1$ and $t>1$.
EXP-$(s,t)$-WT holds for arbitrary $\gamma_j$'s  iff
\begin{equation}\label{eq:expltgt}
\lim_{j\to\infty}\frac{\left(\ln
\wt{\lambda}_j^{-1}\right)^{\eta}} {\ln j}=\infty \ \ \mbox{with}\
\ \eta=\frac{s(t-1)}{t-s}.
\end{equation}

$\bullet$ Let $s<1$ and $t\le 1$.  EXP-$(s,t)$-WT  holds iff
\begin{equation}\label{eq:expq1}
 \lim_{\vz \to 0} \frac{d(\vz)^{1-s} \ln j(\vz)}{(\ln
\vz^{-1})^s} = 0;
\end{equation}

\begin{rem}\label{rem:lambda}
Since $\ln \wt{\lambda}_j^{-1}
=\ln \lambda_j^{-1} + \ln \lambda_1$ and
$\ln \lambda_j^{-1} \to \infty$,
every condition above involving $\{\wt{\lambda}_j\}$
can equivalently be stated with the original sequence
$\{\lambda_j\}$: the additive constant $\ln \lambda_1$ disappears
after division by $\ln j$ or after taking one more logarithm.
\end{rem}

\subsection{Main result}\

Novak and Wo\'zniakowski posed the following open problem
\cite[Open Problem 27]{NW08}: characterize, for weighted linear
tensor product problems under ABS, the various notions of
algebraic tractability in terms of the weights and the univariate
eigenvalues.

We remark that  in the case $\lambda_1 \le 1$  the problem is
settled.  In this paper we solve it completely for
dimension-independent product weights in the case $\lz_1>1$.

\vskip 2mm

\begin{thm}\label{thm:main}
We consider the weighted linear tensor problem $S_{\bsgamma}$  in
the worst case setting for ABS. Assume that  $\lambda_1 > 1$.
Then the following holds. \vspace {2mm}

(i) $S_{\bsgamma}$ is  {\it never} ALG-SPT and {\it never} EXP-SPT
 for arbitrary weight $\gamma_j$'s.\vspace {2mm}

(ii) $S_{\bsgamma}$ is ALG-PT iff $S_{\bsgamma}$ is EXP-PT iff
$\lim\limits_{j \to \infty} \gamma_j = 0$ and
\begin{equation}\label{2.16}
 \limsup_{\vz \to 0}
\frac{d(\vz) \ln j(\vz)}{\ln \ln \vz^{-1}}<\infty.
\end{equation}

(iii) $S_{\bsgamma}$ is ALG-QPT iff $S_{\bsgamma}$ is EXP-QPT iff
$S_{\bsgamma}$ is ALG-PT.\vspace {2mm}

(iv) $S_{\bsgamma}$ is ALG-UWT iff $S_{\bsgamma}$ is EXP-UWT iff
\begin{equation*}
\lim_{j \to \infty} \frac{\ln \ln \lambda_j^{-1}}{\ln \ln j} =
\infty \quad \mbox{and} \quad \lim_{j \to \infty} \frac{\ln \ln
\gamma_j^{-1}}{\ln j} = \infty.
\end{equation*}

(v) $S_{\bsgamma}$ is ALG-$(s,t)$-WT with $s,t>0$ iff
$S_{\bsgamma}$ is EXP-$(t,t)$-WT iff $S_{\bsgamma}$ is
EXP-$(s_1,t)$-WT with $s_1>t$ iff

\quad $\bullet$ for  $t>1$,
\begin{equation*}
\lim_{j \to \infty} \frac{(\ln \lambda_j^{-1})^t}{\ln j} = \infty,
\end{equation*}
with no requirement on the weights;

\quad $\bullet$ for $t=1$,
\begin{equation*}
\lim_{j \to \infty} \gamma_j = 0 \quad \mbox{and} \quad \lim_{j
\to \infty} \frac{\ln \lambda_j^{-1}}{\ln j} = \infty;
\end{equation*}

\quad $\bullet$ for $0<t<1$,
\begin{equation} \label{2.17}
 \quad \lim_{\vz \to 0} \frac{d(\vz)^{1-t} \ln
j(\vz)}{(\ln \vz^{-1})^t} = 0.
\end{equation}

(vi) $S_{\bsgamma}$ is EXP-$(s,t)$-WT with $0<s<t\le1$ iff
\begin{equation*}  \lim_{\vz \to 0}
\frac{d(\vz)^{1-s} \ln j(\vz)}{(\ln \vz^{-1})^s} = 0.
\end{equation*}

(vii) $S_{\bsgamma}$ is EXP-$(s,t)$-WT with $1\le s<t$ iff
$\gamma_j$'s are arbitrary and
$$\lim\limits_{j\to\infty}\frac{(\ln\lambda_j^{-1})^{s}}{\ln
j}=\infty.$$

(viii) $S_{\bsgamma}$ is EXP-$(s,t)$-WT with $0<s<1<t$ for
arbitrary $\gamma_j$'s iff
$$\lim\limits_{j\to\infty}\frac{(\ln\lambda_j^{-1})^{\eta}}{\ln
j}=\infty\ \ {\rm with}\ \ \eta=\frac{s(t-1)}{t-s}.$$

\end{thm}
\vskip 2mm

The conditions \eqref{2.16} and \eqref{2.17} depend on the
``inverse functions'' $j(\vz)$ and $d(\vz)$ of the two sequences
$\{\lz_j\}$ and $\{\gamma_j\}$. Here we give the equivalent
conditions which are described by $\{\lz_j\}$ and $\{\gamma_j\}$
explicitly. See  the following Theorem.

\begin{thm}\label{cor:M}
(i) \eqref{2.16} holds iff
\begin{equation}\label{2.18}
\limsup_{j+k \to \infty} \frac{k \ln j}{\ln \ln
(\gamma_k^{-1}+\lambda_j^{-1})} <\infty.
\end{equation}

(ii) For fixed $t\in(0,1)$, \eqref{2.17} holds iff
\begin{equation}\label{2.19}
  \lim_{j+k \to \infty} \frac{k^{1-t} \ln
j}{\big(\ln (\gamma_k^{-1}+\lambda_j^{-1})\big)^t}= 0.
\end{equation}

\end{thm}

\begin{rem} \label{rem2} In \eqref{2.18} and \eqref{2.19} the sum
$\gamma_k^{-1}+\lambda_j^{-1}$ in the denominator can  be replaced
by the maximum $\max(\gamma_k^{-1},\lz_1\lambda_j^{-1})$, since
\begin{equation}\label{eq:maxsum}
\lz_1^{-1}\max(\gamma_k^{-1},\lz_1\lambda_j^{-1})\le\gamma_k^{-1}+
\lambda_j^{-1}\le2\max(\gamma_k^{-1},\lz_1\lambda_j^{-1}).
\end{equation}
Furthermore,  \eqref{2.18} holds iff there exists a constant $C>0$
for which
\begin{equation*}
 \min(\gamma_k, \lambda_j) \le
\exp(-j^{C k}) , \quad {\rm for\ all}\ \  j,k \in \NN .
\end{equation*}
\end{rem}

\vskip 2mm

The conditions \eqref{2.18} and \eqref{2.19}, though exact, mix
the sequences $\{\lambda_j\}$ and $\{\gamma_j\}$ and may be
difficult to verify. The following corollary provides practical
  necessary conditions and  sufficient conditions on the
weight decay and on the spectral decay.

\begin{cor}\label{cor:explicit}
(i) If  $S_{\bsgamma}$ is ALG-QPT, then
\begin{equation*} \limsup_{j \to \infty} \frac{j}{\ln \ln
\gamma_j^{-1}}<\infty \quad (\mbox{i.e., there exists a } q>1\
{\rm such \ that }\ \gamma_j=O(e^{-q^j})),
\end{equation*}and
\begin{equation*}
\limsup_{j \to \infty} \frac{\ln j}{\ln \ln \lambda_j^{-1}}<\infty
 \quad (\mbox{i.e., there exists a } r>0\ {\rm such \ that }\
\lambda_j=O(e^{-j^r})).
\end{equation*}
Conversely, if there exists a $\theta \in (0,1)$ such that
\begin{equation*}
\limsup_{j \to \infty} \frac{j}{(\ln \ln \gamma_j^{-1})^{\theta}}
<\infty
\quad \mbox{and} \quad
\limsup_{j \to \infty} \frac{\ln j}{(\ln \ln
\lambda_j^{-1})^{1-\theta}}<\infty,
\end{equation*}
then $S_{\bsgamma}$ is EXP-PT.

(ii) Suppose $0<t<1, \,s>0$. If $S_{\bsgamma}$ is ALG-$(s,t)$-WT,
then
\begin{equation*} \lim_{j \to \infty} \frac{j^{1-t}}{(\ln
\gamma_j^{-1})^t}=0 \quad \mbox{and} \quad \lim_{j \to \infty}
\frac{\ln j}{(\ln \lambda_j^{-1})^t}=0.
\end{equation*}
Conversely, if there exists a $\theta \in (0,1)$ such that
\begin{equation*}
\limsup_{j \to \infty} \frac{j^{1-t}}{(\ln
\gamma_j^{-1})^{\theta t}}<\infty \quad \mbox{and}
\quad
\lim_{j \to \infty} \frac{\ln j}{(\ln
\lambda_j^{-1})^{(1-\theta) t}}=0,
\end{equation*}
then $S_{\bsgamma}$ is EXP-$(t,t)$-WT.
\end{cor}

\vskip 2mm

The following example illustrates  Theorem \ref{cor:M} and
Corollary~\ref{cor:explicit}.

\begin{exam}\label{ex:kpw}
Let $\alpha, \beta > 0$ and $\lambda_1 > 1$. Suppose that
\begin{equation*}
\gamma_k = \exp(-\exp(k^{\beta})), \quad \lambda_j = \lambda_1
\exp(-j^{(\ln j)^{\alpha}}).
\end{equation*}
Then
\begin{equation*}
\ln \ln \gamma_k^{-1} = k^{\beta}, \quad \ln \ln \lambda_j^{-1} =
\ln(j^{(\ln j)^{\alpha}}-\ln \lambda_1) = (1+o(1)) (\ln
j)^{\az+1}.
\end{equation*}

If  $\beta < 1$, then $\limsup\limits_{j \to \infty} \frac{j}{\ln
\ln \gamma_j^{-1}}=\infty$ and hence $S_{\bsgamma}$ is not
ALG-QPT.

If $\beta=1$, the necessary condition for ALG-QPT fails. However,
since
\begin{align*}
\limsup_{j+k \to \infty} \frac{k \ln j}{\ln \ln
(\gamma_k^{-1}+\lambda_j^{-1})}&\ge \limsup_{j+k \to \infty}
\frac{k \ln j}{\max(k, (\ln j)^{\az+1})}\\ &\ge \limsup_{k \to
\infty}\frac{k \ln k}{\max(k, (\ln k)^{\az+1})}= \infty,
\end{align*}by Theorem \ref{cor:M} $S_{\bsgamma}$ is not
ALG-QPT.

If $\beta>1$, $\az\ge \frac1{\beta-1}$, then
$$\limsup\limits_{k \to \infty} \frac{k}{(\ln \ln
\gamma_j^{-1})^{\theta}} = 1$$ and
\begin{equation*}
 \limsup_{j \to \infty} \frac{\ln j}{((1+o(1))
(\ln j)^{1+\alpha})^{1-\theta}} = \limsup_{j \to \infty} \frac{\ln
j}{(\ln j)^{(1+\alpha)(1-1/\beta)}}  < \infty,
\end{equation*}where  $\theta = 1/\beta \in (0, 1)$, and hence
$(1+\alpha)(1-1/\beta)\ge1$.

When $0 < \beta \le 1$,  $S_{\bsgamma}$ is not ALG-QPT, but is
EXP-UWT. Since
\begin{equation*}
\frac{\ln \ln \lambda_j^{-1}}{\ln \ln j} = (1+o(1)) \frac{(\ln
j)^{1+\alpha}}{\ln j} \to \infty, \quad \frac{\ln \ln
\gamma_j^{-1}}{\ln j} = \frac{j^{\beta}}{\ln j} \to \infty,
\end{equation*}
\thmref{thm:main} (iv) shows that $S_{\bsgamma}$ is ALG-UWT for
all $\alpha, \beta > 0$.
\end{exam}

\section{Proofs of Theorems \ref{thm:main} and \ref{cor:M} }

In order to prove Theorem \ref{thm:main}, we introduce a series of
 lemmas. The first one converts the ABS problem into the same
problem in additional dimensions; the other ones relate the EXP
tractability notions of $S_{\bsgamma}$ and $\wt{S}_{\bsgamma}$ or
the ALG and EXP tractability  notions of $S_{\bsgamma}$.

\begin{lem}\label{lem:padding}
For all $d, k \in \NN$ and $\vz > 0$, we have
\begin{equation*}
n(\vz, d) \le n(\lambda_1^{k/2} \vz, d+k).
\end{equation*}
\end{lem}

\begin{proof}
For $\bsj\in \NN^d$, we set $$\bsj':=(j_1, \dots, j_d, 1, \dots,
1)\in \NN^{d+k}.$$ It follows from  \eqref{eq:eig2} that $$\uu
(\bsj') = \uu (\bsj), \ \ \gamma_{\uu(\bsj)}=\gamma_{\uu(\bsj')},\
\ {\rm and}\ \ \lambda_{d+k,\bsgamma,\bsj'} =
\lambda_{d,\bsgamma,\bsj} \cdot \lambda_1^k.$$ Therefore,
\begin{equation*}
\begin{aligned}
n(\vz,d)
&=\#\{\bsj \in \NN^d:\lambda_{d,\bsgamma,\bsj}>\vz^2\}\\
&=\#\{\bsj'=(\bsj,1,...,1) \in
\NN^{d+k}:\lambda_{d+k,\bsgamma,\bsj'}>\lambda_1^k \cdot \vz^2\}\\
&\le\#\{\wt{\bsj} \in
\NN^{d+k}:\lambda_{d+k,\bsgamma,\wt{\bsj}}>(\lambda_1^{k/2}
\vz)^2\}\\
&=n(\lambda_1^{k/2} \vz,d+k).
\end{aligned}
\end{equation*}

Lemma \ref{lem:padding} is proved.
\end{proof}

\begin{lem}\label{prop:pt-exppt}
$S_{\bsgamma}$ is EXP-PT iff $\wt{S}_{\bsgamma}$ is EXP-PT.
\end{lem}

\begin{proof}
Assume first that $S_{\bsgamma}$ is EXP-PT. By
\eqref{eq:conv}, we have
\begin{equation*}
\wt{n}(\vz,d)=n(\vz \lambda_1^{d/2},d)\le n(\vz,d),
\end{equation*}
so $\wt{S}_{\bsgamma}$ is EXP-PT as well.

Conversely, assume that $\wt{S}_{\bsgamma}$ is EXP-PT. Then there
exist positive constants $C, p,q$ such that for any $\vz'\in
(0,1)$, we have
$$\wt{n}(\vz', d) \le C d^q (1 + \ln \vz'^{-1})^p.
$$
For $\vz \in (0,1)$ we take $\vz' := \vz \lambda_1^{-d/2} \in
(0,1)$. Equation \eqref{eq:conv} gives,
\begin{align*}
 n(\vz, d) &=  \wt{n}(\vz', d) \le Cd^q (1 +
\ln \vz^{-1} + \frac{d}{2} \ln \lambda_1)^p\\ & \le C_1d^{p+q}(1 +
\ln \vz^{-1})^p.
\end{align*}where $C_1=C(1 + \ln \lambda_1)^p$, and in the last inequality we used the inequality
\begin{equation*}
1 + \ln \vz^{-1} + \frac{d}{2} \ln \lambda_1 \le (1 + \ln
\vz^{-1})(1 + \frac{d}{2} \ln \lambda_1) \le (1 + \ln \lambda_1) d
(1 + \ln \vz^{-1}).
\end{equation*}
This  means $S_{\bsgamma}$ is EXP-PT.

Lemma \ref{prop:pt-exppt} is proved.
\end{proof}

\begin{lem}\label{prop:qpt-pt}
If $S_{\bsgamma}$ is ALG-QPT, then $S_{\bsgamma}$ is EXP-PT.
\end{lem}

\begin{proof}
Assume that there exist constants $t, C$ such that
\begin{equation*}
  n(\vz,d)\le C\exp(t(1+\ln \vz^{-1})(1+\ln d)).
\end{equation*}
 Without
loss of generality, it suffices to consider $\vz \in
(0,\lambda_1^{-1/2})$. Let $K$ be the smallest positive integer
$k$ such that
\begin{equation*}
\lambda_1^{k/2} \vz\ge\lambda_1^{-1/2}, \ \ {\rm equivalently\ }\
k\ge \frac{2 \ln \vz^{-1}}{\ln \lambda_1}-1.
\end{equation*}
Then $$K=\left\lceil \frac{2 \ln \vz^{-1}}{\ln
\lambda_1}\right\rceil-1.$$ We set $\vz'=\lambda_1^{K/2}\vz,\
d'=d+K$. Then $\lambda_1^{-1/2}\le\vz'<1 $, $$1+\ln(\vz')^{-1}\le
1+\frac{\ln\lz_1}2,$$ and \begin{align*}1+\ln d'&\le 1+\ln
\Big(d+\frac{2 \ln \vz^{-1}}{\ln \lambda_1}\Big)\\ &\le\ln \Big(e
(1+\frac2{\ln\lz_1})d(1+\ln \vz^{-1})\Big)\\&=\ln (Ad(1+\ln
\vz^{-1})), \end{align*}where $A=e (1+\frac2{\ln\lz_1})$.

By \lemref{lem:padding}, we have
\begin{equation*}
n(\vz,d)=n(\lambda_1^{-K/2} \vz',d)\le n(\vz',d+K)=n(\vz',d').
\end{equation*}
It follows that
\begin{equation*}
\begin{aligned}
 n(\vz,d)\le n(\vz',d') &\le C\exp(t(1+\ln
(\vz')^{-1})(1+\ln d'))\\
&\le C\exp(t(1+\frac{\ln\lz_1}2) \ln (Ad(1+\ln \vz^{-1}))
)\\
&\le  C'd^p(1+\ln \vz^{-1})^p,
\end{aligned}
\end{equation*}
where $p=t(1+\frac{\ln \lambda_1}2)$ and $C':=C A^p$. This  means
that  $S_{\bsgamma}$ is EXP-PT.

Lemma \ref{prop:qpt-pt} is proved.
\end{proof}

\begin{lem}\label{lem:stwt1}
For fixed $s, t > 0$, $S_{\bsgamma}$ is ALG-$(s,t)$-WT iff
  $S_{\bsgamma}$ is EXP-$(t,t)$-WT.
\end{lem}
\begin{proof}For $s,t>0$, assume that $S_{\bsgamma}$ is
ALG-$(s,t)$-WT. Without loss of generality, it suffices to
consider $\vz \in (0,\lambda_1^{-1/2})$. As in the proof of Lemma
\ref{prop:qpt-pt}, set  $$K=\left\lceil \frac{2 \ln \vz^{-1}}{\ln
\lambda_1}\right\rceil-1,\ \ \vz'=\lambda_1^{K/2}\vz,\ \ d'=d+K.$$
Then $\lambda_1^{-1/2}\le\vz'<1 $ and $$(d')^t=(d+K)^t\le
\Big(d+\frac{2 \ln \vz^{-1}}{\ln \lambda_1}\Big)^t\le
2^{t-1}\Big(d^t+\Big(\frac{2\ln \vz^{-1}} {\ln
\lambda_1}\Big)^t\Big) \le B(d^t+(1+\ln\vz^{-1})^t),
$$where $B=2^{t-1}(1+(2/\ln \lz_1)^t)$. Obviously,
$(\vz')^{-1}+d'\to\infty$ as $\vz^{-1}+d\to \infty$. Since
$S_{\bsgamma}$ is ALG-$(s,t)$-WT, we get
$$\lim_{\vz^{-1}+d\to\infty}\frac{\ln {n}(\vz',d')}{(d')^t+{(\vz')}^{-s}}=0.   $$
 By \lemref{lem:padding}, we have
\begin{equation*}
n(\vz,d)=n(\lambda_1^{-K/2} \vz',d)\le n(\vz',d+K)=n(\vz',d').
\end{equation*}
It follows that \begin{align*}\frac{\ln {n}(\vz,d)}{d^t+(1+\ln
\vz^{-1})^t}&\le \frac{\ln
{n}(\vz',d')}{(d')^t+{(\vz')}^{-s}}\cdot \frac
{(d')^t+{(\vz')}^{-s}}{d^t+(1+\ln \vz^{-1})^t}\\ &\le \frac{\ln
{n}(\vz',d')}{(d')^t+{(\vz')}^{-s}} \cdot \frac
{B(d^t+(1+\ln\vz^{-1})^t)+\lz_1^{s/2}}{d^t+(1+\ln \vz^{-1})^t}.
\end{align*}
Since the second factor in the last term is bounded,  it follows
that
$$\lim_{\vz^{-1}+d\to\infty}\frac{\ln {n}(\vz,d)}{d^t+(1+\ln
\vz^{-1})^t}=0,$$which means that $S_{\bsgamma}$ is
EXP-$(t,t)$-WT.

On the other hand, for fixed $s,t>0$ assume that $S_{\bsgamma}$ is
EXP-$(t,t)$-WT. We have $$\lim_{\vz^{-1}+d\to\infty}\frac{\ln
{n}(\vz,d)}{d^t+(1+\ln \vz^{-1})^t}=0.$$Note that
$$ \frac{\ln
{n}(\vz,d)}{d^t+\vz^{-s}}=\frac{\ln {n}(\vz,d)}{d^t+(1+\ln
\vz^{-1})^t}\cdot \frac{d^t+(1+\ln \vz^{-1})^t}{d^t+\vz^{-s}},$$
and the second factor in the right hand is bounded. It follows
that
$$\lim_{\vz^{-1}+d\to\infty}\frac{\ln
{n}(\vz,d)}{d^t+\vz^{-s}}=0,$$which means that $S_{\bsgamma}$ is
ALG-$(s,t)$-WT.

Lemma \ref{lem:stwt1} is proved.
\end{proof}

\begin{rem}We remark that ALG-$(s,t)$-WT of $S_{\bsgamma}$ depends
only on the parameter $t$ and is independent of the parameter
$s$.\end{rem}

\begin{lem}\label{lem:stwt}
For fixed $s, t > 0$, $S_{\bsgamma}$ is EXP-$(s,t)$-WT iff
 $\wt{S}_{\bsgamma}$ is EXP-$(\min(s,t),t)$-WT.
\end{lem}

\begin{proof}  We first consider the case  $s\le t$.  In this case,
$\min(s,t)=s$. Assume that $S_{\bsgamma}$ is EXP-$(s,t)$-WT. Since
$S_{\bsgamma}$ is harder than $\wt{S}_{\bsgamma}$. It follows that
$\wt{S}_{\bsgamma}$  is also EXP-$(s,t)$-WT.

On the other hand, assume that  $\wt{S}_{\bsgamma}$ is
EXP-$(s,t)$-WT. It follows that
$$\lim_{\vz^{-1}+d\to\infty}\frac{\ln \wt{n}(\vz,d)}{d^t+(1+\ln\vz^{-1})^{s}}
= 0.$$ By \eqref{eq:conv} we have
$${n}(\vz,d)=\wt{n}(\vz\lz_1^{-d/2}, d).
$$It follows that
 $$\frac{\ln
{n}(\vz,d)}{d^t+(1+\ln\vz^{-1})^{s}}=\frac{\ln
\wt{n}(\vz\lz_1^{-d/2},d)}{d^t+(1+\ln(\vz^{-1}\lz_1^{d/2}))^{s}}\cdot
\frac
{d^t+(1+\ln\vz^{-1}+\frac{d}2\ln\lz_1)^{s}}{d^t+(1+\ln\vz^{-1})^{s}}.
$$The first factor in the right hand tends to 0, and the second
factor is bounded. It follows that $S_{\bsgamma}$ is
EXP-$(s,t)$-WT.

Next we consider the case $s>t$. In this case, $\min(s,t)=t$.
Assume that $S_{\bsgamma}$ is EXP-$(s,t)$-WT. It follows that
$S_{\bsgamma}$ is ALG-$(s,t)$-WT. By Lemma \ref{lem:stwt1} we
obtain that $S_{\bsgamma}$ is EXP-$(t,t)$-WT. Hence
$\wt{S}_{\bsgamma}$ is EXP-$(t,t)$-WT.

On the other hand, assume that $\wt{S}_{\bsgamma}$ is
EXP-$(t,t)$-WT. It follows from the proof of the case $s\le t$
that ${S}_{\bsgamma}$ is EXP-$(t,t)$-WT. Hence, ${S}_{\bsgamma}$
is EXP-$(s,t)$-WT for $s>t$.

Lemma \ref{lem:stwt} is proved.
\end{proof}

\begin{rem}From Lemmas \ref{lem:stwt1} and \ref{lem:stwt}, we knows for fixed $s, t>0$, $S_{\bsgamma}$ is
ALG-$(s,t)$-WT  iff ${S}_{\bsgamma}$ is EXP-$(t,t)$-WT iff
${S}_{\bsgamma}$ is EXP-$(s_1,t)$-WT with $s_1> t$  iff
$\wt{S}_{\bsgamma}$ is EXP-$(t,t)$-WT.
\end{rem}

\begin{lem}\label{lem:uwt}\ $S_{\bsgamma}$ is ALG-UWT iff
$S_{\bsgamma}$ is EXP-UWT iff $\wt{S}_{\bsgamma}$ is EXP-UWT.
\end{lem}

\begin{proof} We note that
that $S_{\bsgamma}$ is ALG-UWT iff  $S_{\bsgamma}$ is
ALG-$(t,t)$-WT for all $t\in (0,1)$. By Lemmas \ref{lem:stwt1} and
\ref{lem:stwt}, we obtain that $S_{\bsgamma}$ is ALG-$(t,t)$-WT
for all $t\in (0,1)$ iff $S_{\bsgamma}$ is EXP-$(t,t)$-WT for all
$t\in (0,1)$ iff $\wt{S}_{\bsgamma}$ is EXP-$(t,t)$-WT for all
$t\in (0,1)$. This means that $S_{\bsgamma}$ is ALG-UWT iff
$S_{\bsgamma}$ is EXP-UWT iff $\wt{S}_{\bsgamma}$ is EXP-UWT.

Lemma \ref{lem:uwt} is proved.
\end{proof}

Now we prove \thmref{thm:main}.\vspace{2mm}

\begin{proof}[Proof of \thmref{thm:main}]\

(i). Since EXP-SPT implies ALG-SPT, it suffices to show that
$S_{\bsgamma}$ is not ALG-SPT. Assume that  $S_{\bsgamma}$ is
ALG-SPT. Then there exist two constants $C>0$ and $p\ge 0$ such
that for all $\vz\in (0,1)$ and $d\in\mathbb N$, $$n(\vz,d)\le
C\vz^{-p}.$$ Then for fixed  $\vz_0\in (0,1)$, the sequence
$\{n(\vz_0,d)\}_{d\in\mathbb N}$ is bounded. On the other hand,
consider the multi-indices with a single active coordinate $k\
(1\le k\le d)$ and $j_k = 2$, i.e.,
$\bsj=(j_1,...,j_{k-1},j_k,j_{k+1},...,j_d):=(1,...,1,2,1...,1)$,
then
\begin{equation*}
\lambda_{d,\bsgamma,\bsj}=\gamma_k \lambda_2 \lambda_1^{d-1} >
\vz_0^2 \quad \iff \quad \gamma_k>\vz_0^2
\lambda_1^{1-d}/\lambda_2.
\end{equation*}
The right-hand side tends to $0$ as $d \to \infty$.  Since for any
 $N\in \mathbb N$, $\gamma_N>0$,  we get that there exists a $D\in \mathbb N$ such that $D>N$ and
$$\gamma_N>\vz_0^2 \lambda_1^{1-d}/\lambda_2$$ for all $d>D$. This
means that
$$n(\vz_0,d)= \#\{ \bsj \in \NN^d : \lambda_{d,\bsgamma,\bsj}  >
\vz_0^2 \}\ge\#\{1\le k\le d : \gamma_k>\vz_0^2
\lambda_1^{1-d}/\lambda_2\}\ge N$$ as $d>D$. Since $N$ is
arbitrary, $n(\vz_0,d)$ tends to $\infty$ as $d\to\infty$. This
leads to a contradiction. Hence $S_{\bsgamma}$ is not ALG-SPT.
\vspace{2mm}

(ii), (iii). By Lemma \ref{prop:qpt-pt}, $S_{\bsgamma}$ is ALG-QPT
iff $S_{\bsgamma}$ is EXP-PT.  From the definitions, EXP-PT
implies both ALG-PT and EXP-QPT, and each of the latter in turn
implies ALG-QPT. This means that $S_{\bsgamma}$ is ALG-PT iff
$S_{\bsgamma}$ is EXP-PT iff $S_{\bsgamma}$ is ALG-QPT iff
$S_{\bsgamma}$ is EXP-QPT. Hence (iii) holds. By Lemma
\ref{prop:pt-exppt}, $S_{\bsgamma}$ is EXP-PT iff
$\wt{S}_{\bsgamma}$ is EXP-PT. By \eqref{eq:exppt}, we obtain that
$S_{\bsgamma}$ is EXP-PT iff \begin{equation*} \lim\limits_{j \to
\infty} \gamma_j = 0\ \ {\rm and} \quad \limsup_{\vz \to 0}
\frac{d(\vz) \ln j(\vz)}{\ln \ln \vz^{-1}} < \infty.
\end{equation*}
Hence (ii) holds.\vspace{2mm}

(iv). By \lemref{lem:uwt}, $S_{\bsgamma}$ is ALG-UWT iff
$\wt{S}_{\bsgamma}$ is EXP-UWT. By \eqref{eq:expuwt} we obtain
that $S_{\bsgamma}$ is ALG-UWT iff
\begin{equation*}
\lim_{j \to \infty} \frac{\ln \ln \lambda_j^{-1}}{\ln \ln j} =
\infty \quad \mbox{and} \quad \lim_{j \to \infty} \frac{\ln \ln
\gamma_j^{-1}}{\ln j} = \infty.
\end{equation*}

(v).  By Lemmas \ref{lem:stwt1} and \ref{lem:stwt}, for $s,t>0$,
$S_{\bsgamma}$ is ALG-$(s,t)$-WT  iff $S_{\bsgamma}$ is
EXP-$(t,t)$-WT  iff  $\wt{S}_{\bsgamma}$ is EXP-$(t,t)$-WT. By
\eqref{eq:expsgt}, \eqref{1.8},  and \eqref{eq:expq1}, we get (v).

(vi), (vii), (viii). By Lemma \ref{lem:stwt}, for $0<s<t$,
$S_{\bsgamma}$ is EXP-$(s,t)$-WT   iff  $\wt{S}_{\bsgamma}$ is
EXP-$(s,t)$-WT. By \eqref{eq:expq1}, \eqref{eq:exp1gt} and
\eqref{eq:expsgt}, and \eqref{eq:expltgt} we obtain (vi), (vii),
(viii) respectively.

Theorem \ref{thm:main} is proved.
\end{proof}

\begin{proof}[Proof of Theorem \ref{cor:M}.] \

We remark that any one of the conditions \eqref{2.16},
\eqref{2.17}, \eqref{2.18}, and \eqref{2.19} implies that
$\lim\limits_{j\to\infty}\gamma_j=0$. We write
\begin{equation*}
a_k := \ln \gamma_k^{-1}, \quad b_j := \ln (\lambda_1
\lambda_j^{-1}).
\end{equation*}

(i). Assume that \eqref{2.16} holds. Then there exists a constant
$B>0$ such that \begin{equation*} d(\vz) \ln j(\vz) \le B \ln \ln
\vz^{-1}, \quad \mbox{for all sufficiently small } \vz.
\end{equation*}
Since
$\lim\limits_{j\to\infty}\gamma_j=\lim\limits_{j\to\infty}\lz_j=0$,
we obtain that $d(\vz), j(\vz) \to \infty$ as $\vz \to 0$. By the
definitions
 of  $j(\vz)$ and $d(\vz)$, we have
\begin{equation*}
\gamma_k > \vz^2 \iff k \le d(\vz), \quad {\rm and}\  \ \lambda_j
> \lambda_1 \vz^2 \iff j \le j(\vz).
\end{equation*}
It follows that  $a_{d(\vz)} < 2 \ln \vz^{-1}$ and $b_{j(\vz)} < 2
\ln \vz^{-1}$.

 For every
sufficiently large $j + k$, take $\vz := \vz_{j,k}$ with
\begin{equation*}
\vz_{j,k}^2 := \frac{1}{2} e^{- \max (a_k, b_j)}.
\end{equation*}
This gives
\begin{equation*}
\vz^2 \le \frac{1}{2} e^{- a_k} = \gamma_k/2 < \gamma_k \quad
\mbox{and} \quad \lambda_1 \vz^2 \le \frac{\lambda_1}2 e^{- b_j} =
\frac{\lambda_j}{2} < \lambda_j,
\end{equation*}
It follows that
\begin{equation}\label{eq:dj}
d(\vz) \ge k \quad \mbox{and} \quad j(\vz) \ge j.
\end{equation}
Moreover $\vz \to 0$ as $j + k \to \infty$, and
\begin{equation}\label{eq:lnln}
\ln \ln \vz^{-1} = \ln(\frac{\max(a_k, b_j) +\ln2}2) = (1 + o(1))
\ln \max(a_k, b_j).
\end{equation}
Inserting \eqref{eq:dj} and \eqref{eq:lnln} into \eqref{2.16}
gives
\begin{equation*}
k \ln j \le d(\vz) \ln j(\vz) \le B \ln \ln \vz^{-1} = (B + o(1))
\ln \max(a_k, b_j),
\end{equation*}
that is,
\begin{equation*}
\frac{k \ln j}{\ln \max(a_k, b_j)} \le B + o(1).
\end{equation*}
By Remark \ref{rem2}, \eqref{2.18} follows.

On the other hand, assume that \eqref{2.18} holds.  Then there
exists a constant $B'>0$ such that
\begin{equation}\label{eq:P1}
k \ln j \le B' \ln \max(a_k, b_j) \quad \mbox{for all sufficiently
large } j + k.
\end{equation}
For all sufficiently small $\vz$ we take
\begin{equation*}
k = d(\vz)  \quad \mbox{and} \quad j = j(\vz);
\end{equation*}
then $j + k \to \infty$ as $\vz \to 0$, and
\begin{equation*}
 a_k < 2 \ln \vz^{-1}\  \mbox{ and } \ b_j <
2 \ln \vz^{-1}.
\end{equation*}
It follows  that
\begin{equation*}
\ln \max(a_k, b_j) < \ln(2 \ln \vz^{-1}) = \ln \ln \vz^{-1} + \ln
2.
\end{equation*}
Applying \eqref{eq:P1} to this pair $(k, j)$ gives
\begin{equation*}
d(\vz) \ln j(\vz) = k \ln j \le B' \ln \max(a_k, b_j) \le B' (\ln
\ln \vz^{-1} + \ln 2) \le 2 B' \ln \ln \vz^{-1}
\end{equation*}
for all sufficiently small $\vz$, which is \eqref{2.16}.
\vspace{2mm}

(ii). Assume that \eqref{2.17} holds. Then
\begin{equation*}
d(\vz)^{1-t} \ln j(\vz) = o((\ln \vz^{-1})^t) \quad \mbox{as }\
\vz \to 0.
\end{equation*}
For every sufficiently large $j + k$, take $\vz := \vz_{j,k}$ with
\begin{equation*}
\vz_{j,k}^2:= \frac{1}{2} e^{-\max (a_k, b_j)}.
\end{equation*} as in the proof of the above (i). It follows that
\begin{equation*}k\le d(\vz), \ \ j\le j(\vz), \ \ {\rm and}\ \
\ln \vz^{-1} = \frac{\max(a_k, b_j) +\ln2}2.
\end{equation*}
Inserting into \eqref{2.17} gives
\begin{equation*}
k^{1-t} \ln j \le d(\vz)^{1-t} \ln j(\vz) = o((\ln \vz^{-1})^t) =
o(\max(a_k, b_j)^t),
\end{equation*}
that is,
\begin{equation*}
\frac{k^{1-t} \ln j}{\max(a_k, b_j)^t} \to 0 \quad \mbox{as}\ \ j
+ k \to \infty,
\end{equation*}
which is \eqref{2.19}.

On the other hand, assume that \eqref{2.19} holds. For all
sufficiently small $\vz$ we take
\begin{equation*}
k = d(\vz)  \quad \mbox{and} \quad j = j(\vz);
\end{equation*}
then $j + k \to \infty$ as $\vz \to 0$, and
\begin{equation*}
 a_k < 2 \ln \vz^{-1}\  \mbox{ and } \ b_j <
2 \ln \vz^{-1}.
\end{equation*}
It follows  that
\begin{equation*}
( \max(a_k, b_j))^t< 2^t\,( \ln \vz^{-1})^t .
\end{equation*}
Applying \eqref{2.19} to this pair $(k, j)$ gives
\begin{equation*}
\frac{d(\vz)^{1-t} \ln j(\vz)}{(\ln \vz^{-1})^t} \le \frac{2^t
k^{1-t} \ln j}{\max(a_k, b_j)^t} \to 0 \ \ {\rm as}\ \ \vz\to 0,
\end{equation*}
which is \eqref{2.17}.

Theorem \ref{cor:M} is proved.
\end{proof}

{\small

\noindent\textbf{Funding}

\noindent The authors  were supported by the National Natural
Science Foundation of China (Project no. 12371098).

\

\noindent\textbf{Author contributions}

\noindent All authors contributed equally to this work.

\

\noindent\textbf{Competing interests}

\noindent The authors declare no completing interesting.}

\end{document}